\input amstex
\documentstyle{amsppt}
\magnification 1200

\def\diam{\operatorname{diam}}
\def\R{\Bbb R}

\def\N{\Bbb N}
\def\Z{\Bbb Z}

\def\id{\operatorname{id}}
\def\A{\Cal A}
\def\U{\Cal U}
\def\V{\Cal V}
\def\Ord{\operatorname{Ord}}
\def\cl{\operatorname{cl}}
\def\dim{\operatorname{dim}}
\def\ind{\operatorname{ind}}
\def\Ind{\operatorname{Ind}}
\def\asdim{\operatorname{asdim}}
\def\asInd{\operatorname{asInd}}
\def\asind{\operatorname{asind}}
\def\trInd{\operatorname{trInd}}
\def\trasInd{\operatorname{trasInd}}
\def\trasind{\operatorname{trasind}}
\def\trasdim{\operatorname{trasdim}}
\NoBlackBoxes

\topmatter
\title
On transfinite extension of asymptotic dimension
\endtitle

\author
T.Radul
\endauthor
\address
Dept. de Matematicas, Facultad de Cs. Fisicas y Mat., Universidad
de Concepcion, CASILLA 160-C, Concepcion, Chile
 \newline
e-mail: tarasradul\@yahoo.co.uk
\endaddress

\abstract We prove that a transfinite extension of asymptotic
dimension $\asind$ is trivial.  We introduce a transfinite
extension of asymptotic dimension $\asdim$ and give an example of
metric proper space which has transfinite infinite dimension.
\endabstract
\endtopmatter
\document
\baselineskip18pt

{\bf 0.} Asymptotic dimension $\asdim$ of a metric space was
defined by Gromov for studying asymptotic invariants of discrete
groups [1]. This dimension can be considered as asymptotic
analogue of the Lebesgue covering dimension $\dim$. Dranishnikov
has introduced dimensions $\asInd$ and $\asind$ which are
analogous to large inductive dimension $\Ind$ and small inductive
dimension $\ind$ [2,3].It is known that $\asdim X=\asInd X$ for
each proper metric space with $\asdim X<\infty$. The problem of
coincidence of $\asdim$ and $\asInd$ is still open in the general
case [3].

Extending codomain of $\Ind$ to ordinal numbers we obtain the
transfinite extension $\trInd$ of the dimension $\Ind$. It is
known that there exists a space $S_\alpha$ such that $\trInd
S_\alpha=\alpha$ for each countable ordinal number $\alpha$ [4].
Zarichnyi has proposed to consider transfinite extension of
$\asInd$ and conjectured that this extension is trivial. It is
proved in [5] that if a space has a transfinite asymptotic
dimension $\trasInd$, then this dimension is finite.

We investigate in this paper transfinite extensions for the
asymptotic dimensions $\asind$ and $\asdim$. It appears that
extending codomain of $\asind$ to ordinal numbers we obtain the
trivial extension as well. However, the main result of this paper
is construction of transfinite extension $\trasdim$ of $\asdim$
which is not trivial. Moreover, $\trasdim$ classifies the metric
spaces with asymptotic property $C$ introduced by Dranishnikov
[8].

The paper is organized as follows: in Section 1 we give some
necessary definitions and introduce some denotations, in Section 2
we prove that the transfinite extension of $\asind$ is trivial and
in Section 3 we define the transfinite extension $\trasdim$ of
$\asdim$ and build a proper metric space $X$ such that $\trasdim
X=\omega$.

{\bf 1.} Let $A_1,A_2\subset X$ be two disjoint closed subsets in
a topological space $X$. We recall that a {\it partition} between
$A_1$ and $A_2$ is a subset $C\subset X$ such that there are open
disjoint sets $U_1,U_2$ satisfying the conditions: $X\setminus
C=U_1\cup U_2$, $A_1\subset U_1$ and $A_2\subset U_2$. Clearly a
partition $C$ is a closed subset of $X$.

We will define the asymptotic dimensions $\asind$ and $\asInd$ for
the class of proper metric space. We recall that a metric space is
{\it proper} if every closed ball is compact. We assume that some
base point $x_0\in X$ is chosen for each proper metric space $X$.
The generic metric we denote by $d$. If $X$ is a metric space and
$A\subset X$ we denote by  $B_r(A)$ the open $r$-neighborhood:
$B_r(A)=\{x\in X\mid d(x,A)< r\}$. We call two subsets $A_1,
A_2\subset X$ in a metric space $X$ {\it asymptotically disjoint}
if $\lim_{r\to\infty}d(A_1\setminus B_r(x_0),A_2\setminus
B_r(x_0))=\infty$.

A map $\phi:X\to I=[0,1]$ is called {\it slowly oscillating} if
for any $r>0$, for given $\varepsilon>0$ there exists $D>0$ such
that $\diam \phi(B_r(x))<\varepsilon$ for any $x$ with
$d(x,x_0)\ge D$. If $C_h(X)$ is the set of all continuous slow
oscillating functions $\phi:X\to I$, then the {\it Higson
compactification} is the closure of the image of $X$ under the
embedding $\Phi:X\to I^{C_h(X)}$ defined as
$\Phi(x)=(\phi(x)\mid\phi\in C_h(X))\in I^{C_h(X)}$. We denote the
Higson compactification of a proper metric space $X$ by $cX$ and
the remainder $cX\setminus X$ by $\nu X$. The compactum $\nu X$ is
called {\it Higson corona}. Let us remark that $\nu X$ does not
need to be metrizable.

Let $C$ be a subset of a proper metric space $X$. By $C'$ we
denote the intersection $\cl C\cap \nu X$ of the closure $\cl C$
in the Higson compactification $cX$. Clearly, two sets $A_1$ and
$A_2$ are asymptotically disjoint iff their traces $A_1'$ and
$A_2'$ in the Higson corona are disjoint. Note that for each $r>0$
we have $B_r(C)'=C'$.

Let $A_1, A_2\subset X$ be two asymptotically disjoint subsets of
a proper metric space $X$. A  subset $C\subset X$ is called an
{\it asymptotic separator} for $A_1$ and $A_2$ if its trace $C'$
is a partition for $A_1'$ and $A_2'$ in $\nu X$.

We define $\asInd X=-1$ if and only if $X$ is bounded; $\asInd
X\le n$ if for every two asymptotically disjoint sets $A$,
$B\subset X$ there is an asymptotic separator $C$ with $\asInd
C\le n-1$. Naturally we say $\asInd X=n$ if $\asInd X\le n$ and it
is not true that $\asInd X\le n-1$. We set $\asInd X=\infty$ if
$\asInd X>n$ for each $n\in \N$ [2].

Let  $x\in\nu X$ and $ A\subset X$ such that $x\notin A'$. A
subset $C\subset X$ is called an {\it asymptotic separator} for
$x$ and $A$ if its trace $C'$ is a partition for $\{x\}$ and $A'$
in $\nu X$.

We define $\asind X=-1$ if and only if $X$ is bounded; $\asind
X\le n$ if for every $x\in\nu X$ and $ A\subset X$ such that
$x\notin A'$ there is an asymptotic separator $C$ with $\asind
C\le n-1$. Naturally we say $\asind X=n$ if $\asind X\le n$ and it
is not true that $\asind X\le n-1$. We set $\asind X=\infty$ if
$\asind X>n$ for each $n\in \N$ [3].

There are proved subspace and addition theorems for $\asInd$ in
[5]:

\proclaim {Theorem A} Let $X$ be a proper metric space and
$Y\subset X$. Then $\asInd Y\le\asInd X$.
\endproclaim

\proclaim {Theorem B} Let $X$ be a proper metric space and
$X=Y\cup Z$ where $Y$ and $Z$ are unbounded sets. Then $\asInd
X\le\asInd Y+\asInd Z$.
\endproclaim

Define the transfinite extension $\trasInd X$: $\trasInd X=-1$ if
and only if $X$ is bounded; $\trasInd X\le \alpha$ where $\alpha$
is an ordinal number if for every two asymptotically disjoint sets
$A$, $B\subset X$ there is an asymptotic separator $C$ with
$\trasInd C\le \beta$ for some $\beta<\alpha$. Naturally we say
$\trasInd X=\alpha$ if $\trasInd X\le \alpha$ and it is not true
that $\trasInd X\le \beta$ for some $\beta<\alpha$. We set
$\trasInd X=\infty$ if for each ordinal number $\alpha$ it is not
true that $\trasInd X\le \alpha$. It is proved in [5] that this
extension is trivial:

\proclaim {Theorem C} Let $X$ be a proper metric space such that
$\trasInd X<\infty$. Then $\asInd X<\infty$.
\endproclaim

{\bf 2.} We consider in this section a transfinite extension of
asymptotic dimension $\asind$ and show that this extension is
trivial.

Define the transfinite extension $\trasind X$: $\trasind X=-1$ if
and only if $X$ is bounded; $\trasind X\le \alpha$ where $\alpha$
is an ordinal number if for every $x\in\nu X$ and $ A\subset X$
such that $x\notin A'$ there is an asymptotic separator $C$ with
$\trasind C\le \beta$ for some $\beta<\alpha$. Naturally we say
$\trasind X=\alpha$ if $\trasind X\le \alpha$ and it is not true
that $\trasind X\le \beta$ for some $\beta<\alpha$. We set
$\trasind X=\infty$ if for each ordinal number $\alpha$ it is not
true that $\trasind X\le \alpha$. It follows from the definition
that $\asind X <\infty$ iff $\trasind X<\omega$ where $\omega$ is
the first infinite ordinal number.

\proclaim {Lemma 1} Let $\trasind X=\alpha$ for some ordinal
number $\alpha$. Then for each $\beta<\alpha$ there exists a
subset $X_\beta\subset X$ such that $\trasind X_\beta=\beta$.
\endproclaim

\demo{Proof} We shall apply transfinite induction with respect to
$\alpha$. For $\alpha=0$ the lemma is obvious. Assume that the
theorem holds for all $\alpha<\alpha_0\ge 1$ and consider a proper
metric space $X$ such that $\trasind X=\alpha_0$ as well an
ordinal number $\beta<\alpha_0$. Suppose that $X$ contains no
subset $M$ with $\trasind M=\beta$. By the inductive assumption
$X$ contains no subset $M'$ which satisfies $\beta\le\trasind
M'<\alpha_0$. Thus for every point $x\in \nu X$ and each $A\subset
X$ such that $x\notin A'$ there exists an asymptotic separator $C$
for $x$ and $A$  such that $\trasind C<\beta$. This contradicts,
however, the equality $\trasind X=\alpha_0$, so that $X$ contains
a subset $X_\beta$ with $\trasind X_\beta=\beta$.
\enddemo

It is proved in [3] that $\asind X\le\asInd X$.

\proclaim {Lemma 2} Let $\asind X<\infty$ for some proper metric
space $X$. Then $\asInd X<\infty$ as well.
\endproclaim

\demo{Proof} We use induction with respect to $\asind X$. If
$\asind X=-1$, then $\asInd X=-1$. Suppose we have proved the
lemma for each $i<n\ge 0$. Consider any proper metric space $X$
with $\asind X\le n$. Let $A$ and $B$ be asymptotically disjoint
subsets of $X$ and $a$ is any point of $A'$. Since $\asind X\le
n$, there exists an asymptotic separator $L_a$ between $a$ and $B$
such that $\asind L_a<n$. Consider open disjoint sets $U_a$, $V_a$
in $\nu X$ such that $a\in U_a$, $B'\subset V_a$ and $\nu
X\setminus L_a'=U_a\cup V_a$. Since $A'$ is compact, there exist
points $a_1,\dots,a_k\in A'$ such that $A'\subset\cup_{i=1}^k
U_{a_i}$. Put $U=\cup_{i=1}^k U_{a_i}$, $V=\cap_{i=1}^k V_{a_i}$
and $S=\nu X\setminus (U\cup V)$. Then $U$ and $V$ are open
disjoint subsets of $\nu X$ such that $A'\subset U$ and $B'\subset
V$. Hence $S$ is a partition between $A'$ and $B'$ in $\nu X$.
Moreover, $S=\cup_{i=1}^k L_{a_i}'\setminus U$. Choose a
continuous function $f:\nu X\to [0,1]$ such that
$f(A')\subset\{0\}$ and $f(\nu X\setminus U)\subset\{1\}$. We can
extend this function to a continuous function $F:cX\to [0,1]$. Put
$L=(\cup_{i=1}^k L_{a_i})\setminus(F^{-1}([0,\frac{1}{2}])\cap
X)$. Then we have $S\subset L'$ and hence $L$ is an asymptotic
separator between $A$ and $B$.

Since $\asind L_{a_i}<n$, we have $\asInd L_{a_i}<\infty$ for each
$i$ by inductive assumption. Hence we have $\asInd
L\le\asInd\cup_{i=1}^k L_{a_i}<\infty$ by Theorems A and B. So,
$\trasInd X\le\omega$ and $\asInd X<\infty$ by Theorem C. The
lemma is proved
\enddemo

\proclaim {Theorem 1} Let $\trasind X<\infty$ for some proper
metric space $X$. Then $\asind X<\infty$ as well.
\endproclaim

\demo{Proof} Suppose the contrary. Then there exists a proper
metric space $X$ such that $\trasind X=\alpha$ for some ordinal
number $\alpha\ge\omega$. We can choose a proper metric space $Y$
such that $\trasind Y=\omega$. Let us show that $\asInd Y<\infty$.
Choose any asymptotically disjoint sets $A$ and $B$ in $Y$. Since
$\trasind Y=\omega$, we can choose for each point $a\in A'$ an
asymptotic separator $L_a$ between $a$ and $B$ such that $\asind
L_a<\infty$. So, $\asInd L_a<\infty$ by Lemma 2. Using the same
method as in the proof of Lemma 2, we can choose an asymptotic
separator $L$ between $A$ and $B$ such that $\asInd L<\infty$.
Hence,  $\trasInd Y\le\omega$ and $\asInd Y<\infty$ by Theorem C.
Then $\asind Y\le\asInd Y<\infty$ and we obtain the contradiction.
The theorem is proved.
\enddemo

{\bf 3.} In this section we introduce a transfinite extension of
dimension $\asdim$ introduced by Gromov [1]. A family $\A$ of
subsets of a metric space is called {\it uniformly bounded} if
there exists a number $C>0$ such that $\diam A\le C$ for each
$A\in\A$; $\A$ is called $r$-{\it disjoint} for some $r>0$ if
$d(A_1,A_2)\ge r$ for each $A_1$, $A_2\in\A$ such that $A_1\neq
A_2$.

The {\it asymptotic dimension} of a metric space $X$ does not
exceed $n\in\N\cup\{0\}$ (written $\asdim X\le n$) iff for every
$D>0$ there exists a uniformly bounded cover $\U$ of $X$ such that
$\U=\U_0\cup\dots\cup\U_n$, where all $\U_i$ are $D$-disjoint.
Moreover, we put $\asdim X=-1$ iff $X$ is bounded.

Since the definition of $\asdim$ is not inductive, we cannot
immediately extend this dimension. We need some set-theoretical
construction used by Borst to extend covering dimension and metric
dimension [6,7].

Let $L$ be an arbitrary set. By $Fin L$ we shall denote the
collection of all finite, non-empty subsets of $L$. Let $M$ be a
subset of $Fin L$. For $\sigma\in\{\emptyset\}\cup Fin L$ we put
$$M^\sigma=\{\tau\in Fin L\mid \sigma\cup\tau\in M\ \text{and}\
\sigma\cap\tau=\emptyset\}.$$ Let $M^a$ abbreviate $M^{\{a\}}$ for
$a\in L$.

 Define the ordinal number $\Ord M$ inductively
as follows

$\Ord M=0$ iff $M=\emptyset$,

$\Ord M\le\alpha$ iff for every $a\in L$, $\Ord M^a<\alpha$,

$\Ord M=\alpha$ iff $\Ord M\le\alpha$ and $\Ord M<\alpha$ is not
true, and

$\Ord M=\infty$ iff $\Ord M>\alpha$ for every ordinal number
$\alpha$.

We will need some lemmas from [6]:

\proclaim {Lemma D} Let $L$ be a set and let $M$ be a subset of
$Fin L$. In addition let $n\in\N$. Then $\Ord M\le n$ iff
$|\sigma|\le n$ for each $\sigma\in M$.
\endproclaim

We call a subset $M$ of $Fin L$ \it {inclusive} iff for every
$\sigma$, $\sigma'\in Fin L$ such that $\sigma\in M$ and
$\sigma'\subset \sigma$ also $\sigma'\in M$.

\proclaim {Lemma E} Let $L$ be a set and let $M$ be an inclusive
subset of $Fin L$.  Then $\Ord M=\infty$ iff there exists a
sequence $(a_i)_{i=1}^\infty$ of distinct elements of $L$ such
that $\sigma_n=\{a_i\}_{i=1}^n\in M$ for each $n\in\N$.
\endproclaim

\proclaim {Lemma F} Let $\phi:L\to L'$ be a function and let
$M\subset Fin L$ and $M'\subset Fin L'$ be such that for every
$\sigma\in M$ we have $\phi(\sigma)\in M'$ and
$|\phi(\sigma)|=|\sigma|$. Then $\Ord M\le\Ord M'$.
\endproclaim

Let us define the following collection for a metric space $(X,d)$:
$$A(X,d)=\{\sigma\in Fin \N\mid \ \text{ there is no uniformly bounded  families}\
\V_i\ \text{for}\ i\in\sigma$$ $$\text{such that}\
\cup_{i\in\sigma}\V_i\ \text{covers}\ X\ \text{and}\ \V_i\
\text{is}\ i-\text{disjoint}\}.$$

Let $(X,d)$ be a metric space. Then put $\trasdim X=\Ord A(X,d)$
and \newline $\trasdim X=-1$ iff $X$ is bounded. It follows from
Lemma D that $\trasdim$ is a transfinite extension of $\asdim$:
$\trasdim X\le n$ iff $\asdim X\le n$ for each $n\in\N$.

Dranishnikov has defined asymptotic property $C$ as follows: a
metric space $X$ has asymptotic property $C$ if for any sequence
of natural numbers $n_1<n_2<\dots$ there is a finite sequence of
uniformly bounded families $\{\U_i\}_{i=1}^n$ such that
$\cup_{i=1}^n\U_i$ covers $X$ and $\U_i$ is $n_i$-disjoint [8].

The next proposition follows from Lemma E:

\proclaim {Proposition 1} A metric space $X$ has asymptotic
property $C$ iff $\trasdim X<\infty$.
\endproclaim

\proclaim {Proposition 2} Let $X$ be a metric space and $Y\subset
X$. Then $\trasdim Y\le\trasdim X$.
\endproclaim

\demo{Proof} Put $M=A(Y)$, $M'=A(X)$ and $\phi=\id_\N:\N\to\N$.
Then $M$, $M'$ and $\phi$ satisfy the condition of Lemma F and
$\trasdim Y=\Ord A(Y)\le\Ord A(X)=\trasdim X$.
\enddemo

We are going to construct two examples: a proper metric space
$L_\omega$ such that $\trasdim L_\omega=\omega$ which shows that
this extension is not trivial and  a proper metric space
$L_\infty$ such that $\trasdim L_\infty=\infty$.

We denote by $N_R(A)=\{x\in X|d(x,A)\le R\}$ for a metric space
$X$, $A\subset X$ and $R>0$.

We consider $\Z^n$ with the sup-metric defined as follows
$d((k_1,\dots,k_n),(l_1,\dots,l_n))=\max\{|k_1-l_1|,\dots,|k_n-l_n|\}$.
It follows from [1,Lemma 6.1] that $\asdim \Z^n\le n$. By $k\Z$ we
denote $\{kl|l\in\Z\}$ for some $k\in\N$.

\proclaim {Lemma 3} There exist no uniformly bounded $2k$-disjoint
families $\V_1,\dots,\V_n$ in $(k\Z)^n$ such that
$\cup_{i=1}^n\V_i$ covers $(k\Z)^n$.
\endproclaim

\demo{Proof} Suppose the contrary. Then there exists $s\in\N$ such
that $\diam V\le s$ for each $i\in\{1,\dots,n\}$ and $V\in\V_i$.
We can suppose that $s\ge 2k$. Consider $(k\Z)^n$ as the subset of
$\R^n$. Put $\V=\cup_{i=1}^n\V_i$ and
$\U=\{N_{k/2}(V)\cap[-s,s]^n|V\in\V\}$. Then $\U$ is a finite
closed cover of the cube $[-s,s]^n$ no member of which meets two
opposite faces of $[-s,s]^n$ and each subfamily of $\U$ containing
$n+1$ distinct elements of $\U$ has empty intersection. We obtain
the contradiction with the Lebesgue's Covering Theorem [4, Theorem
1.8.20].
\enddemo

\proclaim {Corollary} $\trasdim(k\Z)^n=\asdim(k\Z)^n=n$ for each
$k,n\in\N$.
\endproclaim

Put $X=\cup_{i=1}^\infty \Z^i$. Define a metric in X. Let
$a=(a_1,\dots,a_l)\in\Z^l$ and $b=(b_1,\dots,b_k)\in\Z^k$. Suppose
that $l\le k$. Consider $a'=(a_1,\dots,a_l,0,\dots,0)\in\Z^k$. Put
$c=0$ if $l=k$ and $c=l+(l+1)+(l+2)+\dots+(k-1)$ if $l<k$. Now,
define $d_\infty(a,b)=\max\{d(a',b),c\}$ where $d$ is sup-metric
in $\Z^k$.

Consider the proper metric space $L_\infty=(X,d_\infty)$ and its
subspace $L_\omega=\cup_{k=1}^\infty(k\Z)^k\subset L_\infty$

\proclaim {Lemma 4} Let $Y$ be a metric space and $X\subset Y$.
Then $\trasdim N_n(X)=\trasdim X$ for each $n\in\N$.
\endproclaim

\demo{Proof} Consider the function $\phi:\N\to\N$ defined as
follows  $\phi(k)=k+2n$ for $k\in\N$. Obviously, we have
$|\phi(\tau)|=|\tau|$ for each $\tau\in Fin\N$. Consider any
$\tau\in\{k_1,\dots,k_l\}\in A(N_n(X))$. Suppose that
$\phi(\tau)=\{k_1+2n,\dots,k_l+2n\}\notin A(X)$. Then there exist
a sequence of uniformly bounded families $\{\U_i\}_{i=1}^l$ such
that $\cup_{i=1}^l\U_i$ covers $X$ and $\U_i$ is $k_i+2n$-disjoint
for each $i\in\{1,\dots,l\}$. Consider the family
$N_n(\U_i)=\{N_n(V)|V\in\U_i\}$ for each $i\in\{1,\dots,l\}$. Then
the families $N_n(\U_1),\dots,N_n(\U_l)$ are uniformly bounded,
$\cup_{i=1}^l N_n(U_i)$ covers $N_n(X)$ and $\U_i$ is
$k_i$-disjoint for each $i\in\{1,\dots,l\}$. We obtain the
contradiction. So, $\phi(\tau)\in A(X)$ and $\trasdim
N_n(X)\le\trasdim X$ by Lemma F. The inequality $\trasdim
N_n(X)\ge\trasdim X$ follows from Proposition 2.
\enddemo

\proclaim {Theorem 2} $\trasdim L_\omega=\omega$.
\endproclaim

\demo{Proof} The inequality $\trasdim L_\omega\ge\omega$ follows
from Proposition 2 and Corollary.

Consider any $n\in\N$. Let us show that $\Ord A(L_\omega)^n\le
n-1$. Consider any $\tau=\{k_1,\dots,k_n\}\in Fin\N$ such that
$n\notin\tau$. It is enough to show that $\tau\cup\{n\}\notin
A(L_\omega)$.

Since $\asdim (n\Z)^n$ and $\cup_{i=1}^n(i\Z)^i\subset
N_{1+2+\dots+n-1}(n\Z^n)$, we have that
$\trasdim\cup_{i=1}^n(i\Z)^i\le n$. Then there exist a sequence of
uniformly bounded families $\{\U_i\}_{i=1}^{n+1}$ such that
$\cup_{i=1}^{n+1}\U_i$ covers $\cup_{i=1}^n(i\Z)^i$, $\U_i$ is
$k_i$-disjoint for each $i\in\{1,\dots,n\}$ and $\U_{n+1}$ is
$n$-disjoint. Consider the family
$\V=\U_{n+1}\cup\{\{x\}|x\in\cup_{i=n+1}^\infty(i\Z)^i\}$. Then
$\V$ is $n$-disjoint uniformly bounded family such that
$(\cup_{i=1}^{n}\U_i)\cup\V$ covers $L_\omega$. Hence
$\tau\cup\{n\}\notin A(L_\omega)$.
\enddemo

\proclaim {Theorem 3} $\trasdim L_\infty=\infty$.
\endproclaim

\demo{Proof} Suppose the contrary. Consider the sequence
$(n_i)_{i=1}^\infty$ where $n_i=i+1$. Then there exists $k\in\N$
such that $\{2,3,\dots,k+1\}\notin A(L_\infty)$. But then
$\{2,3,\dots,k+1\}\notin A(\Z^k)$ and we obtain the contradiction
with Lemma 3.
\enddemo

Finally, we will prove that $\trasdim$ could have only countable
values.

\proclaim {Lemma 5} Let $X$ be a metric space and $\tau\in
Fin\N\cup\{\emptyset\}$ such that $\Ord A(X)^\tau=\alpha$ for some
ordinal number $\alpha$. Then for each $\xi\le\alpha$ there exists
$\sigma\in Fin\N\cup\{\emptyset\}$ such that $\Ord
A(X)^{\tau\cup\sigma}=\xi$.
\endproclaim

\demo{Proof} We shall apply the transfinite induction with respect
to $\alpha$. For $\alpha=0$ the lemma is obvious. Assume that the
lemma holds for all $\alpha<\alpha_0$ and consider a metric space
$X$ and $\tau\in Fin\N\cup\{\emptyset\}$ such that $\Ord
A(X)^\tau=\alpha_0$ as well as an ordinal number $\xi<\alpha_0$.
Suppose that there is no $\sigma\in Fin\N\cup\{\emptyset\}$ such
that $\Ord A(X)^{\tau\cup\sigma}=\xi$. By the inductive assumption
there is no $\sigma'\in Fin\N\cup\{\emptyset\}$ such that
$\xi\le\Ord A(X)^{\tau\cup\sigma'}<\alpha_0$. Then $\Ord
A(X)^{\tau\cup\{n\}}<\xi$ for each $n\in\N\setminus\tau$ and we
obtain the contradiction with $\Ord A(X)^\tau=\alpha_0$.
\enddemo

\proclaim {Theorem 4} If we have  $\trasdim X<\infty$ for some
metric space $X$, then $\trasdim X<\omega_1$.
\endproclaim

\demo{Proof} Suppose the contrary. Then there exists a metric
space $X$such that $\trasdim X\ge\omega_1$. We can choose $\tau\in
Fin\N\cup\{\emptyset\}$ such that $\Ord A(X)^\tau=\omega_1$. Then
for each $n\in\N\setminus\tau$ we have $\Ord
(A(X)^\tau)^n=\xi_n<\omega_1$. Then
$\omega_1=\sup\{\xi_n|n\in\N\setminus\tau\}$ and we obtain the
contradiction.
\enddemo

\enddocument